\documentclass[12pt,a4paper]{article}
\usepackage[a4paper, textheight=24cm,textwidth=16cm]{geometry}

\usepackage{amssymb}
\usepackage{amstext}
\usepackage{amsmath}
\usepackage{amsthm}
\usepackage{amscd}
\usepackage{tikz-cd}
\usepackage{color}
\usepackage{url}
\usepackage{mathrsfs}
\usepackage{hyperref}
\usepackage[OT2,T1]{fontenc}

\DeclareMathOperator{\Spec}{Spec\,}

\DeclareMathOperator{\Frac}{Frac}

\newcommand{\val}[1]{\vert#1\vert}
\newcommand{\Val}[1]{\Vert#1\Vert}

\newcommand{\Ll}{V}

\newcommand{\inj}{\hookrightarrow}

\newcommand{\surj}{\twoheadrightarrow}

\newcommand{\ZZ}{{\mathbb Z}}

\newcommand{\NN}{{\mathbb N}}

\newcommand{\cB}{{\mathscr{B}}}

\newcommand{\cI}{{\mathscr I}}

\def\Fp{\mathfrak{p}}
\def\FP{\mathfrak{P}}

\def\Fm{\mathfrak{m}}
\def\Fn{\mathfrak{n}}

\def\Fp{\mathfrak{p}}
\def\Fq{\mathfrak{q}}

\newcommand{\cO}{{\mathscr{O}}}

\newcommand{\ol}{\overline}

\newcounter{nc}
\renewcommand{\thenc}{{\rm(\roman{nc})}}
\newenvironment{romlist}%
{\begin{list}{\thenc}{
\usecounter{nc}
\parsep=0pt
\setlength  \labelwidth{\leftmargin}
\addtolength\labelwidth{-\labelsep}
}
}{\end{list}}
\newcommand{\pauseromlist}%
{\global\edef\savecount{\arabic{nc}}\end{romlist}}
\newcommand{\finpauseromlist}%
{\begin{romlist}\setcounter{nc}{\savecount}}
\newcounter{nnc}
\renewcommand{\thennc}{{\rm(\alph{nnc})}}
{\begin{list}{\thennc}{
\usecounter{nnc}
\parsep=0pt
\setlength  \labelwidth{\leftmargin}
\addtolength\labelwidth{-\labelsep}
}
}{\end{list}} 
{\begin{list}{$\bullet$}{\parsep=0pt}}{\end{list}} 
\newcommand{\dem}{\noindent{\slshape D\'emonstration. }}

\newcounter{ctnum}
\renewcommand{\thectnum}{\textup{(\arabic{ctnum})}}
\newenvironment{numlist}%
{\begin{list}{\thectnum}{
\usecounter{ctnum} 
\parsep=0pt
\leftmargin=0pt%
\setlength{\itemindent}{\labelwidth}%
\addtolength{\itemindent}{\labelsep}%
}
}{\end{list}}

\usepackage[french, english]{babel}

\newtheorem{teointro}{Th\'eor\`eme}

\newtheorem{MBteo}[subsection]{Th\'eor\`eme}
\newtheorem{MBprop}[subsection]{Proposition}

\newtheorem{MBdefi}[subsection]{D\'efinition}
\theoremstyle{definition}
\newtheorem{MBrem}[subsection]{Remarque}

\newtheorem{MBsubrem}[subsubsection]{Remarque}
\newtheorem{MBsubrems}[subsubsection]{Remarques}

\theoremstyle{plain}

\newtheorem{MBsubprop}[subsubsection]{Proposition}

\newtheorem{MBsubcor}[subsubsection]{Corollaire}

\theoremstyle{remark}

\theoremstyle{plain}

\newcommand\rref[1]{{\rm\ref{#1}}}

\newcommand{\red}{\mathrm{red}}

\DeclareMathOperator{\fdim}{fdim}
\DeclareMathOperator{\Tdim}{Tor-dim}
\newcommand{\moins}{\smallsetminus}

\title{Une construction d'extensions faiblement non ramifi\'ees d'un anneau de valuation}
\author{%
Laurent Moret-Bailly\thanks{Univ Rennes, CNRS, IRMAR - UMR 6625, F-35000 Rennes, France}
 \thanks{\tt laurent.moret-bailly[AT]univ-rennes1.fr}
\thanks{L'auteur a b\'en\'efici\'e du soutien du projet Geolie (ANR-15-CE 40-0012) de l'Agence nationale de la recherche. }
}

\begin{document}
\selectlanguage{french}
\date{21 mai 2022}
\maketitle
%

\begin{center}
\fbox{ 
\parbox{12cm}{\begin{center}Paru dans {\sl Rend. Sem. Mat. Univ. Padova} (\href{https://ems.press/journals/rsmup/articles/5898517}{Online first})\smallskip\\
\url{https://doi.org/10.4171/rsmup/94}
\end{center}}%
}
\end{center}

\selectlanguage{english}
\begin{abstract}
Given a valuation ring $\Ll$, with residue field $F$ and value group $\Gamma$, we give a sufficient condition for a local ring dominating $\Ll$ to be a valuation ring with the same value group. When $\Ll$ contains a field $k$, we apply this result to the problem of constructing a valuation ring $W$ containing $\Ll$ and a prescribed extension $k'$ of $k$, with value group $\Gamma$ and residue field generated by $k'$ and $F$; this is possible in particular when either $k'$ or $F$ is separable over $k$.
\end{abstract}
\selectlanguage{french}
\begin{abstract}
\'Etant donn\'e un anneau de valuation $\Ll$, de corps r\'esiduel $F$ et de groupe des valeurs $\Gamma$, on donne une condition suffisante pour qu'un anneau local dominant $\Ll$ soit un anneau de valuation de groupe $\Gamma$. Lorsque $\Ll$ contient un corps $k$, ce r\'esultat est appliqu\'e \`a la construction d'un anneau de valuation  contenant $\Ll$ et une extension donn\'ee $k'$ de $k$, de groupe $\Gamma$ et de corps r\'esiduel engendr\'e par $k'$ et $F$. Cela s'av\`ere possible, notamment, lorsque $k'$ ou $F$ est s\'eparable sur $k$.
\end{abstract}
\noindent{\sl Classification AMS 2010:} 13F30, 12J10.

\section{Introduction}
Si $\Ll$ est un anneau de valuation, une $\Ll$-alg\`ebre $W$ est une \emph{extension faiblement non ramifi\'ee} de $\Ll$ si $W$ est un anneau de valuation dominant $\Ll$ et de m\^eme groupe des valeurs que $\Ll$. Le principal r\'esultat de ce travail est le suivant (th\'eor\`eme \ref{teo:ppal}):
\begin{teointro}\label{thmA}
Soit $\Ll$ un anneau de valuation, d'id\'eal maximal $\Fm$, et soit $W$ une $\Ll$-alg\`ebre locale et plate, d'id\'eal maximal $\Fn$. On suppose que:
\begin{romlist}
\item\label{thmA1} $\Fn=\Fm W$;
\item\label{thmA2} $W$ est la colimite d'un syst\`eme inductif filtrant $(R_{i})_{i\in I}$ de $\Ll$-alg\`ebres plates essentiellement de type fini, dont les morphismes de transition $u_{i,j}:R_{i}\to R_{j}$ $(i\leq j)$ sont plats.
\end{romlist} 
Alors $W$ est une extension faiblement non ramifi\'ee de $\Ll$.
\end{teointro}
La d\'emonstration figure au \S\ref{sec:thppal}. On se ram\`ene facilement au cas o\`u $W=A_{\Fp}$ o\`u $A$ est une $\Ll$-alg\`ebre plate de type fini et o\`u $\Fp$ est un point maximal de $\Spec(A/\Fm A)$. Par localisation \'etale, gr\^ace \`a un r\'esultat de Gruson et Raynaud, on peut m\^eme supposer que $A$ est un $\Ll$-module libre et que $\Fp=\Fm A$. Dans ce cas, on construit explicitement la valuation associ\'ee \`a $W$ \`a partir de la \og norme\fg\ naturelle sur le $\Ll$-module libre $A$.
\medskip

Le th\'eor\`eme \ref{thmA} est ensuite appliqu\'e au probl\`eme suivant, qui a motiv\'e ce travail (\`a partir de questions li\'ees au crit\`ere valuatif de propret\'e et aux vari\'et\'es \og pseudo-compl\`etes\fg, voir \cite{LMBCompact}): on part d'un anneau de valuation $\Ll$, de corps r\'esiduel $F$ et contenant un corps $k$. On donne de plus une extension $k'$ de $k$. Il s'agit de construire un anneau de valuation dominant $\Ll$ et contenant $k'$, dont le corps r\'esiduel et le groupe des valeurs soient \og aussi proches que possible\fg\ de ceux de $\Ll$, compte tenu de ces contraintes. Le r\'esultat le plus utile est alors le suivant (th\'eor\`eme \ref{th:ExtCorpsRed}, dont l'\'enonc\'e est plus pr\'ecis):
\begin{teointro}\label{thmB}
Avec les notations pr\'ec\'edentes, on suppose de plus que l'anneau $k'\otimes_{k}F$ est \emph{r\'eduit}. Alors  il existe une extension faiblement non ramifi\'ee $W$ de $\Ll$ contenant $k'$ et  dont le corps r\'esiduel $F_{1}$ est une extension compos\'ee de $k'$ et $F$ sur $k$; de plus on peut imposer que $F_{1}$ soit s\'eparable sur $F$ (resp.\ sur $k'$) si $k'$ (resp.\ $F$) est s\'eparable sur $k$.
\end{teointro}

Pour construire $W$, on part de l'anneau $\Ll':=k'\otimes_{k}\Ll$, et l'on choisit un id\'eal premier minimal $\ol{\Fp}$ de $\Ll'/\Fm\Ll'=k'\otimes_{k}F$. Si $\Fp$ est l'id\'eal premier de $\Ll'$ correspondant, on pose $W=\Ll'_{\Fp}$. Il est alors facile de v\'erifier que le th\'eor\`eme \ref{thmA} s'applique \`a $W$, la condition \ref{thmA1} utilisant l'hypoth\`ese que $k'\otimes_{k}F$ est r\'eduit.

\begin{MBrem}
Les th\'eor\`emes \ref{thmA} et \ref{thmB} sont faciles lorsque $V$ est un anneau de valuation discr\`ete: on utilise le fait que tout anneau local noeth\'erien dont l'id\'eal maximal est en\-gen\-dr\'e par un \'el\'ement r\'egulier est un anneau de valuation discr\`ete. La preuve de ce cas particulier du th\'eor\`eme \ref{thmB} est contenue dans celle de \cite[Proposition C.1.1]{CGP}.
\end{MBrem}
Sans supposer   $k'\otimes_{k}F$  r\'eduit, on obtient un peu moins (th\'eor\`eme \ref{th:ExtCorpsGen}, plus pr\'ecis):

\begin{teointro}\label{thmC}
Les notations sont celles du th\'eor\`eme \rref{thmB}, sans l'hypoth\`ese sur $k'\otimes_{k}F$; on note $p$ l'exposant caract\'eristique de $k$ et $\Gamma$ le groupe de $\Ll$.

Alors il existe un anneau de valuation $W$ dominant $\Ll$ et contenant $k'$, avec les propri\'et\'es suivantes, o\`u $\Delta$ d\'esigne le groupe de $W$ et $F_{1}$ son corps r\'esiduel:
\begin{romlist}
\item\label{thmC1}  le groupe $\Delta/\Gamma$ est de $p$-torsion;
\item\label{thmC2} $F_{1}$ est radiciel sur son sous-corps engendr\'e par $k'$ et $F$.
\end{romlist}
\end{teointro}

\noindent{\sl Plan de l'article.} Le \S\ref{sec:GenVal} donne des rappels sur les valuations. Au \S\ref{sec:AlgLibres} on \'etudie la norme naturelle sur un $\Ll$-module libre (dans toute la suite, $\Ll$ d\'esignera un anneau de valuation) puis sur une $\Ll$-alg\`ebre libre comme $\Ll$-module; le r\'esultat final (proposition \ref{prop:AlgLibreIntegre}), bien que tout \`a fait \'el\'ementaire, est essentiel pour la preuve du th\'eor\`eme principal. Ce dernier est \'etabli au \S\ref{sec:thppal}, o\`u figure aussi une autre d\'emonstration, indiqu\'ee par H. Knaf, du fait que (sous les hypoth\`eses du th\'eor\`eme \ref{thmA}) $W$ est un anneau de valuation. Le \S\ref{sec:ExtComp} traite des propri\'et\'es de certaines extensions compos\'ees dites \og (strictement) ma\-xi\-males\fg, qui apparaissent naturellement dans la preuve du th\'eor\`eme \ref{thmB}. Au \S\ref{sec:CorpsCtes} on d\'emontre les th\'eor\`emes \ref{th:ExtCorpsRed} et \ref{th:ExtCorpsGen}, dont les th\'eor\`emes \ref{thmB} et \ref{thmC} sont des sous-produits.
\medskip

\noindent{\sl Conventions.} Les anneaux sont commutatifs et unitaires. Si $x$ est un point d'un sch\'ema $X$, on note $\cO_{X,x}$ son anneau local et $\kappa(x)$ son corps r\'esiduel.\medskip

\noindent{\sl Remerciements.} L'auteur remercie Hagen Knaf pour ses commentaires et en particulier pour les \'enonc\'es \ref{cor:CorpsResSep} et \ref{prop:Homol}.

\section{Valuations: g\'en\'eralit\'es}\label{sec:GenVal}
\subsection{Notations et conventions}\label{ssec:NotVal}
On utilise la notation multiplicative pour les valuations. Si $(\Gamma,\leq)$ est un groupe ab\'elien (multiplicatif) totalement ordonn\'e, on note $\ol{\Gamma}$ le mono\"{\i}de $\Gamma\cup\{0\}$ o\`u $0.\gamma=0$ et $0\leq\gamma$ pour tout $\gamma\in\ol{\Gamma}$. Une valuation $v$ sur un corps $K$, \`a valeurs dans $\Gamma$, est un morphisme de mono\"{\i}des  $v:(K,\times)\to\ol{\Gamma}$, qui sera plut\^ot not\'e $z\mapsto\val{z}_{v}$ ou simplement $z\mapsto\val{z}$, v\'erifiant $\val{0}=0$, $\val{K^\times}\subset\Gamma$, et $\val{z+t}\leq\max\{\val{z},\val{t}\}$ pour tous $z$ et $t$ dans $K$. L'\emph{anneau} de $v$ est la \og boule unit\'e\fg\ $\Ll=\{z\in K\colon\val{z}\leq1\}$ et son \emph{groupe} est l'image de $K^\times$ par $v$, qui s'identifie au quotient $K^\times/\Ll^\times$; en g\'en\'eral, ce groupe sera suppos\'e \'egal \`a $\Gamma$.

Pour les propri\'et\'es g\'en\'erales des valuations, nous renvoyons \`a \cite{BourAC5-7}, chapitre 6 (qui utilise la notation additive).

\subsection{Extensions faiblement non ramifi\'ees} 
Si $\Ll$ est un anneau de valuation, une \emph{extension} de $\Ll$ est un morphisme  $j:\Ll\to W$ o\`u $W$ est un anneau de valuation et o\`u $j$ est injectif et local. Une telle extension induit un morphisme injectif et croissant de groupes to\-ta\-le\-ment ordonn\'es $\val{j}:\Gamma\inj\Delta$ entre les groupes associ\'es. Nous dirons que $j$ est \emph{faiblement non ramifi\'e} (terminologie de \cite[Tag 0ASG]{StProj}) si $\val{j}$ est un isomorphisme. 

\begin{MBsubprop}\label{prop:FaibNonRam} Soit $\Ll$ un anneau de valuation.
\begin{numlist}
\item\label{prop:FaibNonRam0}  Le hens\'elis\'e de $\Ll$ est une extension faiblement non ramifi\'ee de $\Ll$.
\item\label{prop:FaibNonRam1} Si $j:\Ll\to W$ est une extension faiblement non ramifi\'ee de $\Ll$, l'application naturelle $\FP\mapsto j^{-1}(\FP)$ est un hom\'eomorphisme de $\Spec(W)$ sur  $\Spec(\Ll)$; son inverse est donn\'e par $\Fp\mapsto\Fp W$.
\item\label{prop:FaibNonRam2} Soit $W$ une $\Ll$-alg\`ebre qui est r\'eunion filtrante d'extensions faiblement non ramifi\'ees de $\Ll$. Alors $\Ll\to W$ est une extension faiblement non ramifi\'ee.
\end{numlist}
\end{MBsubprop} 
\dem (D\'etails laiss\'es au lecteur) Pour \ref{prop:FaibNonRam0}, voir \cite[Tag 0ASK]{StProj} ou \cite[Remark 6.1.12\,(vi)]{GabRam}. 
La partie \ref{prop:FaibNonRam1} d\'ecoule facilement du dictionnaire entre id\'eaux d'un anneau de valuation et sous-groupes du groupe associ\'e, et la partie \ref{prop:FaibNonRam2}  est imm\'ediate.\qed

\section{Alg\`ebres libres sur un anneau de valuation}\label{sec:AlgLibres}
On fixe un anneau de valuation $\Ll$, de corps des fractions $K$, d'id\'eal maximal $\Fm$, de corps r\'esiduel $F=\Ll/\Fm$, de groupe $\Gamma=K^\times/\Ll^\times$; conform\'ement aux conventions du \S\ref{sec:GenVal}, on notera $\val{z}\in\ol{\Gamma}$ la valuation de $z\in K$.

\subsection{Norme sur un $\Ll$-module libre}\label{ssec:NormeLibre}
Soit $E$ un $\Ll$-module libre, vu comme un sous-$\Ll$-module de  $E_{K}:=K\otimes_{\Ll}E$. Si $\mathscr{B}=(e_{i})_{i\in I}$ est une $\Ll$-base de $E$  (donc une $K$-base de $E_{K}$), et si $z\in E_{K}$ s'\'ecrit $z=\sum_{i\in I}z_{i}e_{i}$, o\`u les $z_{i}$ sont dans $K$, on pose $\Val{z}_{E}:=\max\limits_{i\in I}\val{z_{i}}\in\ol{\Gamma}$. 

\begin{MBsubprop}\label{prop:GenLibre}
\begin{numlist}
\item\label{prop:GenLibre1} Pour tout $z\in E_{K}$, on a $\Val{z}_{E}=\min\{\val{\alpha}\colon\alpha\in K, z\in\alpha E\}$. En particulier $\Val{z}_{E}$ est ind\'ependant du choix de $\cB$.
\item\label{prop:GenLibre2} Pour $z$, $z'$ dans $E_{K}$ et $\alpha$ dans $K$, on a les propri\'et\'es suivantes:\\
\setcounter{nc}{0}
\noindent\makebox[.495\textwidth][l]{\qquad\refstepcounter{nc}\phantom{(iii)}\llap{\thenc\ }\label{prop:GenLibre21}  $\Val{z+z'}_{E}\leq\max(\Val{z}_{E},\Val{z'}_{E})$;}
\makebox[.495\textwidth][l]{\qquad\refstepcounter{nc}\phantom{(iii)}\llap{\thenc\ }\label{prop:GenLibre22}  $\Val{z}_{E}=0\Leftrightarrow z=0$;}
\makebox[.495\textwidth][l]{\qquad\refstepcounter{nc}\phantom{(iii)}\llap{\thenc\ }\label{prop:GenLibre23}   $\Val{z}_{E}\leq1\Leftrightarrow z\in E$;}
\makebox[.495\textwidth][l]{\qquad\refstepcounter{nc}\phantom{(iii)}\llap{\thenc\ }\label{prop:GenLibre24}  $\Val{z}_{E}<1\Leftrightarrow z\in \Fm E$;}
\makebox[.495\textwidth][l]{\qquad\refstepcounter{nc}\phantom{(iii)}\llap{\thenc\ }\label{prop:GenLibre25}  $\Val{\alpha z}_{E}=\val{\alpha}\,\Val{z}_{E}$.}
\item\label{prop:GenLibre3} Si $E$ n'est pas nul, tout $z\in E_{K}$ peut s'\'ecrire $z=\alpha z_{1}$ avec $\alpha\in K$ et $\Val{z_{1}}_{E}=1$ \emph{(de sorte que $\val{\alpha}=\Val{z}_{E}$)}.
\item\label{prop:GenLibre4} Soit $b:E\times E'\to G$ une application $\Ll$-bilin\'eaire, o\`u $E$, $E'$ et $G$ sont des $\Ll$-modules libres. Pour tout $z\in E$ et tout $z'\in E'$, on a $\Val{b(z,z')}_{G}\leq\Val{z}_{E}\,\Val{z'}_{E'}$.
\end{numlist}
\end{MBsubprop}

La d\'emonstration est laiss\'ee au lecteur (pour \ref{prop:GenLibre4}, on pourra remarquer que pour $\alpha$, $\alpha'\in\Ll$ on a $b(\alpha E\times\alpha' E')\subset \alpha\alpha' G$).\qed

\subsection{Norme: le cas d'une $\Ll$-alg\`ebre}\label{ssec:NormeAlgebreLibre}

Nous allons appliquer ce qui pr\'ec\`ede en prenant pour $E$ une $\Ll$-alg\`ebre $A$, suppos\'ee non nulle et libre comme $\Ll$-module; on identifie $\Ll$ \`a un sous-anneau de $A$ par le morphisme structural. On pose $A_{K}:=K\otimes_{\Ll} A$, $\ol{A}:=F\otimes_{\Ll} A=A/\Fm A$. Pour $z\in A_{K}$, on notera simplement $\Val{z}$ l'\'el\'ement $\Val{z}_{A}\in\ol{\Gamma}$.

\begin{MBsubprop}\label{prop:GenAlgLibre}
\begin{numlist}
\item\label{prop:GenAlgLibre1} Pour tout $\alpha\in K$, on a $\Val{\alpha}=\val{\alpha}$. 
\item\label{prop:GenAlgLibre2} $A\cap K=\Ll$.
\item\label{prop:GenAlgLibre3} Pour tous $z$ et $z'$ dans $A_{K}$, on a $\Val{zz'}\leq\Val{z}\,\Val{z'}$. 
\item\label{prop:GenAlgLibre4} Pour tout $z\in A^\times$, on a $\Val{z}=1$. 
\end{numlist}
\end{MBsubprop}
\dem \ref{prop:GenAlgLibre1} D'abord, comme $A$ est $\Ll$-libre et non nul, on a $\ol{A}\neq0$ et en particulier $1\in A\moins\Fm A$. D'apr\`es \ref{prop:GenLibre}\,\ref{prop:GenLibre2}, \ref{prop:GenLibre23} et \ref{prop:GenLibre24}, ceci implique que $\Val{1}=1$. Pour $\alpha\in K $ quelconque, on a donc, appliquant \ref{prop:GenLibre}\,\ref{prop:GenLibre2}\,\ref{prop:GenLibre25}, $\Val{\alpha}=\Val{\alpha.1}=\val{\alpha}\,\Val{1}=\val{\alpha}$.

L'assertion \ref{prop:GenAlgLibre2} en r\'esulte: si $\alpha\in A\cap K$, on a $\Val{\alpha}\leq1$ donc $\val{\alpha}\leq1$. Enfin, \ref{prop:GenAlgLibre3}  est  un cas particulier de \ref{prop:GenLibre}\,\ref{prop:GenLibre4}, et implique \ref{prop:GenAlgLibre4}  puisque $\Val{1}=1$.\qed

\begin{MBsubcor}\label{cor:AlgLibreRed}
Si $\ol{A}$ est r\'eduit, $A$ est r\'eduit.
\end{MBsubcor}
\dem Soit $z\in A$; \'ecrivons $z=\alpha z_{1}$ avec $\alpha\in \Ll$ et $\Val{z_{1}}=1$. Supposons $z$ nilpotent et non nul. Alors $\alpha\neq0$, donc $z_{1}$ est nilpotent puisque $\alpha$ est r\'egulier dans $A$. Par suite $\ol{z_{1}}$ est nilpotent dans $\ol{A}$, et n'est pas nul puisque $\Val{z_{1}}=1$, donc $\ol{A}$ n'est pas r\'eduit.\qed 
\medskip

Le corollaire  \ref{cor:AlgLibreRed} a la cons\'equence suivante, signal\'ee par Hagen Knaf:

\begin{MBsubcor}\label{cor:CorpsResSep} On suppose que $\Ll$ contient un corps $k$ tel que $F$ soit s\'eparable sur $k$. Alors, pour tout id\'eal premier $\Fp$ de $\Ll$, le corps r\'esiduel $\kappa(\Fp)$ est s\'eparable sur $k$.
\end{MBsubcor}
\dem Rempla\c{c}ant $\Ll$ par $\Ll/\Fp$, on peut supposer que $\Fp=0$: il s'agit donc de montrer que l'extension $K/k$ est s\'eparable ou, de fa\c{c}on \'equivalente, que $k'\otimes_{k}K$ est r\'eduit pour toute extension $k'$ de $k$. Fixons $k'$ et posons $\Ll'=k'\otimes_{k}\Ll$: alors $\Ll'$ est une $\Ll$-alg\`ebre libre comme $\Ll$-module, et $\Ll'/\Fm\Ll'$ s'identifie \`a   $k'\otimes_{k}F$ donc est r\'eduit vu l'hypoth\`ese sur $F$. On d\'eduit alors de \ref{cor:AlgLibreRed} que $\Ll'$ est r\'eduit, donc que son localis\'e $k'\otimes_{k}K$ l'est aussi.\qed
\medskip

La proposition suivante g\'en\'eralise \cite[Lemma 6.1.13]{GabRam} (qui traite le cas o\`u $A$ est un $\Ll$-module de type fini):

\begin{MBsubprop}\label{prop:AlgLibreIntegre}
On suppose que $\ol{A}$ est int\`egre. Alors $A$ est int\`egre, et $\Val{.}$ se prolonge en une valuation sur $\Frac(A)$, prolongeant $\val{.}$, de m\^eme groupe $\Gamma$, d'anneau $A_{\Fm A}$ et de corps r\'esiduel $\Frac(\ol{A})$.
\end{MBsubprop}
\dem Soient $z=\alpha z_{1}$ et $u=\beta u_{1}$ non nuls dans $A$, avec $\alpha$, $\beta\in \Ll$ et $\Val{z_{1}}=\Val{u_{1}}=1$. Alors $\ol{z_{1}}$ et  $\ol{u_{1}}$ ne sont pas nuls dans $\ol{A}$ d'apr\`es \ref{prop:GenLibre}\,\ref{prop:GenLibre2}\,\ref{prop:GenLibre24}, donc $\ol{z_{1}u_{1}}\neq0$. Ceci montre que $zu\neq0$ (donc $A$ est int\`egre), mais aussi que $\Val{z_{1}u_{1}}=1$, de sorte que  $\Val{zu}=\val{\alpha\beta}=\Val{z}\,\Val{u}$. Compte tenu des propri\'et\'es g\'en\'erales vues en \ref{prop:GenLibre}, cela montre que $\Val{.}$ se prolonge en une valuation sur $\Frac(A)$ (encore not\'ee $\Val{.}$), prolongeant $\val{.}$ et \`a valeurs dans $\ol{\Gamma}$.

Posons $L=\Frac(A)=\Frac(A_{K})$ et soit $R=\{x\in L\colon \Val{x}\leq1\}$ l'anneau de  $\Val{.}$. Il est clair que $A\subset R$. Posons $\Fp=\Fm A=\{z\in A\colon\Val{z}<1\}$. Si $z\in A\moins\Fp$, on a $\Val{z}=1$, donc $z\in R^\times$: ceci montre que $A_{\Fp}\subset R$. In\-ver\-se\-ment, soit $x\in R$: \'ecrivons $x={\alpha z}/{u}$ avec $\alpha\in K$, $z\in A$, $u\in A$ et $\Val{z}=\Val{u}=1$ (donc $z,\, u\in A\moins\Fp$). On a alors $\Val{x}=\val{\alpha}\leq1$ donc $\alpha z\in A$, d'o\`u $x\in A_{\Fp}$. Ainsi, $R=A_{\Fp}$, dont le corps r\'esiduel est bien $\Frac(\ol{A})$. \qed

\section{Le th\'eor\`eme principal}\label{sec:thppal}
\begin{MBteo}\label{teo:ppal} Soit $\Ll$ un anneau de valuation, d'id\'eal maximal $\Fm$, et soit $W$ une $\Ll$-alg\`ebre locale d'id\'eal maximal $\Fn$. On suppose que:
\begin{romlist}
\item\label{teo:ppal2} $\Fn=\Fm W$;
\item\label{teo:ppal1} $W$ est colimite d'un syst\`eme inductif filtrant $(W_{i})_{i\in I}$ de $\Ll$-alg\`ebres plates essentiellement de type fini, dont les morphismes de transition $u_{i,j}:W_{i}\to W_{j} $ $(i\leq j)$ sont plats.
\end{romlist} 
Alors $W$ est un anneau de valuation faiblement non ramifi\'e sur $\Ll$.
\end{MBteo}
\dem 

\noindent{\sl a) R\'eduction au cas hens\'elien.}  Le hens\'elis\'e  $\Ll^h$ de $\Ll$ est un anneau de valuation fai\-ble\-ment non ramifi\'e sur $\Ll$ (\ref{prop:FaibNonRam}\,\ref{prop:FaibNonRam0}). Posons $W'=\Ll^h\otimes_{\Ll}W$: alors la $\Ll^h$-alg\`ebre $W'$ v\'erifie encore \ref{teo:ppal2} et \ref{teo:ppal1}, et est une $W$-alg\`ebre locale fid\`element plate; si l'on montre que $W'$ est un anneau de valuation faiblement non ramifi\'e sur $\Ll^h$, la m\^eme propri\'et\'e en r\'esultera pour $\Ll\inj W$ par descente (noter en particulier que l'application $J\mapsto JW'$ est une injection strictement croissante de l'ensemble $\cI$ des id\'eaux de $W$ vers l'ensemble $\cI'$ des id\'eaux de $W'$, de sorte que si $\cI'$ est totalement ordonn\'e par inclusion il en est de m\^eme de $\cI$). Nous supposerons donc dans la suite que $\Ll$ est hens\'elien.
\smallskip

\noindent{\sl b) R\'eduction au cas essentiellement de type fini.} Pour tout $i\in I$, soit $\Fn_{i}\subset W_{i}$ l'image r\'eciproque de $\Fn$: 
rempla\c{c}ant le syst\`eme $(W_{i})$ par $\left((W_{i})_{\Fn_{i}}\right)$, on peut supposer 
que chaque $W_{i}$ est local et que les morphismes de transition $u_{i,j}$ sont locaux (donc aussi les $u_{i}: W_{i}\to W$);  les  $u_{i,j}$ sont alors fid\`element plats et donc injectifs. Par \ref{prop:FaibNonRam}\,\ref{prop:FaibNonRam2}, il suffit de voir que chaque $W_{i}$ est une extension faiblement non ramifi\'ee de $\Ll$. 
Or, par fid\`ele platitude, on a $\Fm W_{i}= \Fm W\cap W_{i}=\Fn_{i}$ pour tout $i$, donc $W_{i}$ v\'erifie l'hypoth\`ese \ref{teo:ppal2} de l'\'enonc\'e; nous sommes donc ramen\'es \`a montrer le cas particulier o\`u $W$ est lui-m\^eme une $\Ll$-alg\`ebre locale, plate et essentiellement de type fini, d'id\'eal maximal $\Fm W$.
\smallskip

\noindent{\sl c) Fin de la d\'emonstration.} On suppose donc que  $W=A_{\Fp}$, o\`u $A$ est une $\Ll$-alg\`ebre de type fini et $\Fp$ un id\'eal premier de $A$. Soit $J\subset A$ le sous-$\Ll$-module de torsion de $A$, qui est un id\'eal de $A$: alors $JA_{\Fp}=0$ puisque $W$ est plat sur $\Ll$, donc $W=(A/J)_{\Fp}$ et l'on peut (rempla\c{c}ant $A$ par $A/J$) supposer $A$ plat sur $V$. Ceci implique que $A$ est une $\Ll$-alg\`ebre \emph{de pr\'esentation finie}, d'apr\`es  \cite[Theorem 3]{Nag66} (ou {\cite[Corollaire 3.4.7]{GruRa71}}).

Posons $\ol{A}=A/\Fm A$ et $\ol{\Fp}:=\Fp/\Fm A=\Fp\ol{A}$: alors  $\ol{A}_{\ol{\Fp}}=W/\Fm W$ est le corps r\'esiduel de $W$, donc $\ol{\Fp}$ est un id\'eal premier minimal de $\ol{A}$, et $\Spec(\ol{A})$ est r\'eduit au point $\ol{\Fp}$; comme $\ol{A}$ est noeth\'erien on peut, en localisant $A$, supposer que $\ol{A}$ est int\`egre et que $\Fp=\Fm A$. 

D'apr\`es \cite[cor. 3.3.13]{GruRa71}, il existe $f\in A\moins\Fp$ tel que $A_{f}$ soit un $\Ll$-module projectif (l'hypoth\`ese hens\'elienne intervient ici). Puisque $\Ll$ est local, $A_{f}$ est m\^eme libre (\cite[Theorem 2]{Kap}, ou \cite[Tag 0593]{StProj}).
Rempla\c{c}ant $A$ par $A_{f}$, nous sommes ramen\'es au cas o\`u $A$ est libre sur $\Ll$, qui r\'esulte de \ref{prop:AlgLibreIntegre}.\qed
\begin{MBsubrem} Si l'on affaiblit l'hypoth\`ese \ref{teo:ppal2} du th\'eor\`eme \ref{teo:ppal} en supposant seulement que $\Fn=\sqrt{\Fm W}$ (autrement dit, $\dim(W/\Fm W)=0$), on peut montrer, en utilisant les r\'esultats de \cite{Nag66}, que les fibres du morphisme ca\-no\-nique $\pi:\Spec(W)\to\Spec(\Ll)$ sont de dimension $0$; de plus, si  $W$ est uni\-branche,  $\pi$ est un hom\'eomorphisme.
\end{MBsubrem}

\subsection{La voie homologique}\label{ssec:Homol}
Les r\'esultats de cette section nous ont \'et\'e indiqu\'es, pour l'essentiel,  par Hagen Knaf. 
Si $R$ est un anneau, nous noterons $\fdim(R)\in\NN\cup\{+\infty\}$ la \emph{dimension faible} de $R$, d\'efinie comme la borne sup\'erieure des tor-dimensions de tous les $R$-modules. On rappelle \cite[Corollary 4.2.6]{Glaz} (ou \cite[Tag092S]{StProj}) qu'un anneau local $R$ est un anneau de valuation si et seulement si $\fdim(R)\leq1$.
Rappelons aussi qu'un anneau $R$ est \emph{coh\'erent} si tout id\'eal de type fini de $R$ est de pr\'esentation finie comme $R$-module.

\begin{MBsubprop}\label{prop:Homol} Soit $\Ll$ un anneau de valuation, d'id\'eal maximal $\Fm$, et soit $W$ une $\Ll$-alg\`ebre locale et plate, d'id\'eal maximal $\Fn$. On suppose que:
\begin{romlist}
\item\label{prop:Homol1} $\Fn=\Fm W$;
\item\label{prop:Homol2} $W$ est un anneau coh\'erent.
\end{romlist} 
Alors $W$ est un anneau de valuation.
\end{MBsubprop}
\dem Montrons que $\fdim(W)\leq1$. On a l'in\'egalit\'e 
\cite[Theorem 3.1.3]{Glaz}
$$\fdim(W)\leq \fdim(W/\Fn)+\Tdim_{W}(W/\Fn)$$
(valable d\`es que $W$ est un anneau coh\'erent et $\Fn$ un id\'eal contenu dans le radical de $W$). 
On a \'evidemment $\fdim(W/\Fn)=0$ puisque $W/\Fn$ est un corps. D'autre part, comme $\Ll$ est un anneau de valuation, on a $\Tdim_{\Ll}(\Ll/\Fm)=1$, d'o\`u  
$$\Tdim_{W}(W/\Fn)=\Tdim_{W}(W/\Fm W)\leq1$$
 puisque $W$ est plat sur $\Ll$, d'o\`u la conclusion.\qed

\begin{MBsubrems}\begin{numlist}\label{rem:coh}
\item Si $\Ll$ est un anneau de valuation, toute $\Ll$-alg\`ebre de pr\'esentation finie (et en particulier toute $\Ll$-alg\`ebre plate de type fini, d'apr\`es \cite[Theorem 3]{Nag66}) est un anneau coh\'erent \cite[Theorem 7.3.3]{Glaz}. D'autre part, si $(R_{i})_{{i\in I}}$ est un syst\`eme inductif filtrant d'anneaux coh\'erents, \`a morphismes de transition plats, l'anneau $\varinjlim_{i}R_{i}$ est encore coh\'erent  \cite[Theorem 2.2.3]{Glaz}.

\hskip1.5em Par suite, dans le  th\'eor\`eme \ref{teo:ppal}, la condition \ref{teo:ppal1} entra\^{\i}ne directement que $W$ est  coh\'erent. La proposition \ref{prop:Homol} red\'emontre donc que $W$ est un anneau de valuation.

\item En revanche, dans la situation de \ref{prop:Homol}, soient $\Gamma$ et $\Delta$ les groupes de valuation respectifs de $\Ll$ et $W$. L'injection $\Gamma\subset\Delta$ induit une injection $\Gamma^{<1}\subset\Delta^{<1}$ entre \'el\'ements $<1$ des deux groupes. La condition \ref{prop:Homol1} \'equivaut \`a dire que $\Gamma^{<1}$ est \emph{cofinal} dans $\Delta^{<1}$, et n'implique pas en g\'en\'eral que $\Gamma=\Delta$; si $\Gamma\subsetneqq\Delta$, on d\'eduit donc du th\'eor\`eme \ref{teo:ppal} que $W$ n'est pas colimite filtrante de $\Ll$-alg\`ebres plates essentiellement de type fini, \`a morphismes de transition plats.
\end{numlist}
\end{MBsubrems}

\section{Compl\'ements sur les extensions compos\'ees}\label{sec:ExtComp}

Cette section est ind\'ependante des pr\'ec\'edentes.
\begin{MBdefi}\label{def:PtStrMax} Soit $X$ un sch\'ema. Un point $x$ de $X$ est dit \emph{strictement maximal} s'il v\'erifie les conditions \'equivalentes suivantes:
\begin{romlist}
\item\label{def:PtStrMax1} l'anneau local $\cO_{X,x}$ est un corps;
\item\label{def:PtStrMax2} $\cO_{X,x}$ est r\'eduit et $x$ est un point maximal de $X$;
\item\label{def:PtStrMax3} le morphisme canonique $\Spec(\kappa(x))\to X$ est plat.
\end{romlist}
\end{MBdefi}

Soient $L$ et $M$ deux extensions d'un corps $K$. Rappelons \cite[V, \S2]{BourA4-7} qu'une \emph{extension compos\'ee} de $L$ et $M$ est un triplet $(K\subset E,u,v)$ o\`u $E$ est une extension de $K$ et o\`u $u:L\to E$, $v:M\to E$ sont deux $K$-plongements dont la r\'eunion des images engendre $E$. On en d\'eduit un morphisme de $K$-alg\`ebres $L\otimes_{K}M\xrightarrow{u\otimes v}E$, et un point $x$ de $X:=\Spec(L\otimes_{K}M)$ \`a corps r\'esiduel isomorphe \`a $E$; le morphisme $u\otimes v$ se factorise en $L\otimes_{K}M\to \cO_{X,x}\surj E$. Cette construction fournit une bijection entre  l'ensemble sous-jacent \`a $X$ et l'ensemble des classes d'isomorphie d'extensions compos\'ees de $L$ et $M$.

\begin{MBdefi}\label{def:maxi} Soit $(E,u,v)$ une extension compos\'ee de $L$ et $M$, et soit $x$ le point de  $X:=\Spec(L\otimes_{K}M)$ correspon\-dant.
Nous dirons que $(E,u,v)$ est \emph{$K$-maximale} (resp. \emph{strictement $K$-maximale}) si $x$  est un point maximal (resp. strictement maximal) de $X$.
\end{MBdefi}
\begin{MBsubrems} \begin{numlist}
\item Nous dirons \og (strictement) maximale\fg\ au lieu de \og (strictement) $K$-maximale\fg\ s'il n'y a pas d'ambigu\"{\i}t\'e sur $K$. Toutefois, on prendra garde que si $E$ est compos\'ee de $L$ et $M$ sur $K$, elle l'est \'egalement sur tout sous-corps de $K$, alors que la maximalit\'e n'est pas conserv\'ee en g\'en\'eral.
\item La maximalit\'e admet une interpr\'etation plus classique, qui ne sera pas utilis\'ee dans la suite: avec les notations de  \ref{def:maxi}, pour que $E$ soit $K$-maximale, il faut et il suffit que $u(L)$ et $v(M)$ soient alg\'ebriquement disjointes sur $K$, c'est-\`a-dire \cite[V,\S14]{BourA4-7} que toute partie de $u(L)$ al\-g\'e\-bri\-que\-ment libre sur $K$ soit encore alg\'ebriquement libre sur $v(M)$.
\end{numlist}
\end{MBsubrems}

\begin{MBprop}\label{prop:ExtComp} Avec $K$, $L$, $M$ et $X$ comme ci-dessus, soient $x$ un point de   $X$ et $(E=\kappa(x),u,v)$ l'extension compos\'ee correspondante.
\begin{numlist}
\item\label{prop:ExtComp0} Si l'extension $L/K$ est radicielle, alors $x$ est le seul point de  $X$,  et $E$ (\'egale \`a $(L\otimes_{K}M)_{\red}$) est, \`a isomorphisme unique pr\`es, la seule extension compos\'ee de $L$ et $M$ sur $K$.
\item\label{prop:ExtComp2} Soit $L_{0}$ un sous-corps de $L$ contenant $K$. Si $E$ est $K$-maximale, il en est de m\^eme du sous-corps $E_{0}\subset E$ engendr\'e par $L_{0}$ et $M$, 
vu comme extension compos\'ee de $L_{0}$ et $M$.
\item\label{prop:ExtComp3} Soit $K_{0}$ un sous-corps de $K$. Si $E$ est $K_{0}$-maximale, alors elle est $K$-maximale; la r\'eciproque est vraie si $K$ est alg\'ebrique sur $K_{0}$.
\item\label{prop:ExtComp5} On suppose que $E$ est maximale. Si l'extension $L/K$ est s\'eparable, alors $E$ est strictement maximale, et est une extension s\'eparable de $M$.
\end{numlist}
\end{MBprop}
\dem L'assertion \ref{prop:ExtComp0}  est bien connue: voir par exemple \cite[V, \S3, exercice 9]{BourA4-7}. 
Montrons \ref{prop:ExtComp2}.
 Soit $x_{0}\in\Spec(L_{0}\otimes_{K}M)$ le point correspondant \`a $E_{0}$. Alors $x_{0}$ est l'image de $x$ par le morphisme naturel $\Spec(L\otimes_{K}M)\to\Spec(L_{0}\otimes_{K}M)$. Ce morphisme est plat, donc envoie un point maximal sur un point maximal, d'o\`u l'assertion.
\smallskip

\noindent\ref{prop:ExtComp3} Soit $z\in X_{0}:=\Spec(L\otimes_{K_{0}}M)$ le point correspondant \`a $E$ vue comme extension compos\'ee sur $K_{0}$. 
On a un diagramme cart\'esien de sch\'emas
$$\begin{tikzcd}
\llap{$x\in$}X=\Spec(L\otimes_{K}M) \dar\rar[hook]{j}
	&X_{0}=\Spec(L\otimes_{K_{0}}M)\rlap{$\ni z$}\dar{f} \\
\Spec(K) \rar[hook] 
	&\Spec(K\otimes_{K_{0}}K)
\end{tikzcd}$$
o\`u $z=j(x)$. Comme les fl\`eches horizontales sont des immersions ferm\'ees, il est clair que si $z$ est maximal dans $X_{0}$, alors $x$ est maximal dans $X$, d'o\`u la premi\`ere assertion. Pour la r\'eciproque, on observe que $f$ est plat; par suite, pour que $z$ soit maximal dans $X_{0}$, il faut et il suffit qu'il soit maximal dans $f^{-1}(f(z))$ (autrement dit, que $x$ soit maximal dans $X$) et que $f(z)$ soit maximal dans $\Spec(K\otimes_{K_{0}}K)$, ce qui est automatique si $K$ est alg\'ebrique sur $K_{0}$ puisqu'alors $\Spec(K\otimes_{K_{0}}K)$ est de dimension $0$.
\smallskip

\noindent\ref{prop:ExtComp5} Si $L/K$ est s\'eparable, alors $L\otimes_{K}M$ est r\'eduit, ainsi que  son localis\'e  $\cO_{X,x}$, de sorte que $E$ est strictement maximale. En outre,  $L\otimes_{K}M$ (et donc aussi  $\cO_{X,x}$) est m\^eme une $M$-alg\`ebre \emph{g\'eom\'etriquement} r\'eduite, d'o\`u l'assertion de s\'eparabilit\'e.\qed

\begin{MBsubrem}
Dans l'assertion \ref{prop:ExtComp5}, pour conclure que $E/M$ est s\'eparable, on ne peut se passer de l'hypoth\`ese que $E$ est maximale. Par exemple, soit $K$ non parfait de caract\'eristique $p$ et soit $a\in K\moins K^p$. Consid\'erons $E=K(a^{1/p},x)$ o\`u $x$ est une ind\'etermin\'ee. Alors $E$ est compos\'ee de $L=K(x)$ et $M=K(x+a^{1/p})$; $L$ est \'evidemment s\'eparable sur $K$, mais $E=M(a^{1/p})$ n'est pas s\'eparable sur $M$.
\end{MBsubrem}

\section{Extension avec sous-corps impos\'e}\label{sec:CorpsCtes}
\subsection{Notations}\label{ssec:NotCorpsCtes}
On garde les notations du \S \ref{sec:AlgLibres}; on suppose en outre que $\Ll$ contient un corps $k$, et l'on fixe une extension $k'$ de $k$. On pose $\Ll'=k'\otimes_{k}\Ll$ et $F'=k'\otimes_{k}F=\Ll'/\Fm\Ll'$. On note $p$ l'exposant caract\'eristique de $k$.

On fixe de plus un id\'eal premier \emph{minimal} $\ol{\Fp}$ de $F'$, et l'on note $\Fp$ l'id\'eal premier de $\Ll'$ correspondant, de sorte que $\ol{\Fp}=\Fp/\Fm\Ll'$. On d\'esigne par $x$ le point de $\Spec(\Ll')$ (ou de $\Spec(F')$) correspondant \`a $\Fp$, et par 
$$W:=\Ll'_{\Fp}$$
son anneau local dans $\Spec(\Ll')$; on a un diagramme cart\'esien de sch\'emas affines
$$\begin{tikzcd}
\makebox[0cm][r]{$x\in$}\Spec(F'_{\ol{\Fp}}) \dar{\ol{j}}\rar[hookrightarrow] 
	& \Spec({W}) \dar{j} & \\
\Spec(F') \dar{\ol{\pi}}\rar[hookrightarrow] 
	& \Spec(\Ll') \dar{\pi} \rar
	& \Spec(k') \dar \\
\Spec(F) \rar[hookrightarrow] 
	& \Spec(\Ll)  \rar
	& \Spec(k)
\end{tikzcd}
$$
o\`u les fl\`eches $j$ et $\ol{j}$ sont les localisations au point $x$, qui est un point maximal de $\Spec(F')$; les fl\`eches du type $\inj$ sont des immersions ferm\'ees. Pour tout point $y$ de $\Spec(\Ll')$, l'ex\-ten\-sion $\kappa(y)/k$ est compos\'ee de $k'$ et de $\kappa(\pi(y))$, puisque la fibre $\pi^{-1}(\pi(y))$ s'identifie \`a $\Spec(k'\otimes_{k}\kappa(\pi(y)))$ (voir les rappels du \S \ref{sec:ExtComp}). 

Il est clair que ${W}$ est une $\Ll$-alg\`ebre locale fid\`element plate; son corps r\'esiduel $\kappa(x)=(F'_{\ol{\Fp}})_{\red}$ est, vu le choix de $x$, une extension compos\'ee \emph{maximale} de $k'$ et $F$. 

\begin{MBteo}\label{th:ExtCorpsRed} \textup{(le cas strictement maximal)} Les hypoth\`eses et notations \'etant celles de \rref{ssec:NotCorpsCtes}, on suppose de plus que $x$ est un point \emph{strictement maximal} \textup{(\ref{def:PtStrMax})} de $\Spec(F')$. Alors:
\begin{numlist}
\item\label{th:ExtCorpsRed1} ${W}$ est une extension faiblement non ramifi\'ee de $\Ll$.
\item\label{th:ExtCorpsRed2} Le morphisme naturel $\pi\circ j:\Spec({W})\to\Spec(\Ll)$ est un hom\'eomorphisme; de plus, pour tout id\'eal premier $\Fq$ de $\Ll$, l'unique id\'eal premier de ${W}$ au-dessus de $\Fq$ est $\Fq':=\Fq{W}$, et son corps r\'esiduel $\kappa(\Fq')$ est une extension compos\'ee \emph{stric\-te\-ment maximale} de $k'$ et $\kappa(\Fq)$ sur $k$.
\item\label{th:ExtCorpsRed3} Si l'extension $k'/k$ est s\'eparable, il en est de m\^eme de $\kappa(x)/F$, et plus g\'en\'eralement (avec les notations de \rref{th:ExtCorpsRed2}) de $\kappa(\Fq')/\kappa(\Fq)$, pour tout $\Fq\in\Spec(\Ll)$.
\item\label{th:ExtCorpsRed4} Si l'extension $F/k$ est s\'eparable, il en est de m\^eme de $\kappa(x)/k'$, et plus g\'en\'eralement de $\kappa(\Fq')/k'$, pour tout $\Fq'\in\Spec({W})$.
\end{numlist} 
\end{MBteo}
\dem  \ref{th:ExtCorpsRed1} Montrons que  les hypoth\`eses du th\'eor\`eme \ref{teo:ppal} sont satisfaites. Par hypoth\`ese, l'anneau local $F'_{\ol{\Fp}}$ co\"{\i}ncide avec $\kappa(\Fp)$; cela \'equivaut \`a dire que l'id\'eal maximal de ${W}$ est engendr\'e par $\Fm$, ce qui est l'hypoth\`ese 
\ref{teo:ppal2} 
de \ref{teo:ppal}. V\'erifions l'hypoth\`ese \ref{teo:ppal1}: soit $(k_{i})_{i\in I}$ la famille ordonn\'ee des sous-extensions de $k'$ qui sont de type fini sur $k$. Alors ${W}$ est limite inductive de la famille $\left((k_{i}\otimes_{k}\Ll)_{\Fp_{i}}\right)_{i\in I}$ o\`u $\Fp_{i}$ est l'id\'eal premier de $k_{i}\otimes_{k}\Ll$ au-dessous de $\Fp$; les morphismes de transition sont plats puisque les $k_{i}$ sont des corps, et chaque $(k_{i}\otimes_{k}\Ll)_{\Fp_{i}}$ est essentiellement de type fini sur $\Ll$ puisque les $k_{i}$ sont des $k$-alg\`ebres essentiellement de type fini. On peut donc appliquer le th\'eor\`eme \ref{teo:ppal}, ce qui donne \ref{th:ExtCorpsRed1}.
\medskip

\noindent\ref{th:ExtCorpsRed2} est cons\'equence de  \ref{th:ExtCorpsRed1} et de la proposition \ref{prop:FaibNonRam}\,\ref{prop:FaibNonRam1}, \`a l'exception des assertions sur les extensions compos\'ees. Pour $\Fq$ et $\Fq'$ comme dans l'\'enonc\'e, on a $\kappa(\Fq')=\kappa(\Fq)\otimes_{\Ll}{W}$ (toujours d'apr\`es  \ref{prop:FaibNonRam}\,\ref{prop:FaibNonRam1}) et en particulier $\kappa(\Fq')$ est plat sur $\kappa(\Fq)\otimes_{\Ll}\Ll'\cong\kappa(\Fq)\otimes_{k}k'$, donc $\Fq'$ est un point strictement maximal de $\Spec(\kappa(\Fq)\otimes_{k}k')$ (condition \ref{def:PtStrMax3} de \ref{def:PtStrMax}), d'o\`u la conclusion.
\medskip

\noindent\ref{th:ExtCorpsRed3} est cons\'equence de  \ref{th:ExtCorpsRed2}, compte tenu de \ref{prop:ExtComp}\,\ref{prop:ExtComp5}. Dans \ref{th:ExtCorpsRed4}, la s\'eparabilit\'e de $\kappa(x)/k'$ r\'esulte aussi de \ref{prop:ExtComp}\,\ref{prop:ExtComp5}, et celle des autres corps r\'esiduels s'en d\'eduit par \ref{cor:CorpsResSep}.\qed

\begin{MBteo}\label{th:ExtCorpsGen} \textup{(le cas g\'en\'eral)} Soient $\Ll$, $k$ et $k'$ comme dans  \rref{ssec:NotCorpsCtes}. Il existe un anneau de valuation $W$ dominant $V$ et contenant $k'$, avec les propri\'et\'es suivantes:
\begin{romlist}
\item\label{th:ExtCorpsGen1} si $\Delta$ d\'esigne le groupe de $W$, le quotient $\Delta/\Gamma$ est de $p$-torsion;
\item\label{th:ExtCorpsGen2} le morphisme naturel $\Spec({W})\to\Spec(\Ll)$ est un hom\'eomorphisme;
\item\label{th:ExtCorpsGen3} pour tout $\Fq'\in\Spec(W)$, d'image $\Fq$ dans $\Spec(\Ll)$, $\kappa(\Fq')$ est radiciel sur son sous-corps engendr\'e par $k'$ et $\kappa(\Fq)$; ce dernier est une extension compos\'ee maximale de $k'$ et $\kappa(\Fq)$ sur $k$.
\end{romlist} 
\end{MBteo}
\dem 
Pour fixer les id\'ees, on supposera que $p>1$ (en caract\'eristique nulle, le th\'eor\`eme \ref{th:ExtCorpsRed} s'applique). Pour tout anneau $R$ de caract\'eristique $p$, on notera $R^\dag=R^{p^{-\infty}}$ sa cl\^oture parfaite. L'anneau $\Ll^\dag$ est un anneau de valuation, de groupe $\Gamma^\dag:=\ZZ[1/p]\otimes_{\ZZ}\Gamma$; explicitement, 
$$\Ll^\dag=\{z\in K^\dag \colon  z^{p^n}\in\Ll\text{ pour $n\in\NN$ assez grand}\}.$$
Soit $k''$ l'unique extension compos\'ee de $k'$ et $k^\dag$ (cf. \ref{prop:ExtComp}\,\ref{prop:ExtComp0}). Comme $k^\dag\subset\Ll^\dag$ et que $k^\dag$ est parfait, on peut appliquer le th\'eor\`eme \ref{th:ExtCorpsRed} \`a $k^\dag\subset\Ll^\dag$ et \`a l'extension $k''/k^\dag$: on obtient un anneau de valuation $W$ contenant $k''$, extension faiblement  non ramifi\'ee de $\Ll^\dag$ et donc de groupe $\Gamma^\dag$. Ceci prouve les assertions  \ref{th:ExtCorpsGen1} et  \ref{th:ExtCorpsGen2}.

Soient $\Fq$ et $\Fq'$ comme dans \ref{th:ExtCorpsGen3}, et soit $\Fq^\dag=\Fq'\cap\Ll^\dag$. Alors $\Fq^\dag$ est le point de $\Spec(\Ll^\dag)$ au-dessus de $\Fq$, de sorte que $\kappa(\Fq^\dag)$ s'identifie \`a $\kappa(\Fq)^\dag$. D'apr\`es \ref{th:ExtCorpsRed}, le corps r\'esiduel $\kappa(\Fq')$ est extension compos\'ee maximale, sur $k^\dag$, de $k''=k'k^\dag$ et  $\kappa(\Fq)^\dag$; il est donc radiciel sur l'extension compos\'ee $k'\kappa(\Fq)$. Il reste \`a voir que celle-ci est $k$-maximale: or on vient de voir que $(\kappa(\Fq'), k'', \kappa(\Fq)^\dag)$ est $k^\dag$- maximale, et donc $k$-maximale d'apr\`es \ref{prop:ExtComp}\,\ref{prop:ExtComp3} puisque $k^\dag$ est alg\'ebrique sur $k$. Comme $k'\subset k''$ et $\kappa(\Fq)\subset\kappa(\Fq)^\dag$, on d\'eduit de  \ref{prop:ExtComp}\,\ref{prop:ExtComp2} que $k' \kappa(\Fq)$ est bien compos\'ee $k$-maximale de $k'$ et $\kappa(\Fq)$.\qed


\end{document}